\documentclass[10pt,twoside]{article}
\usepackage{graphicx}
\usepackage{psfig}
\usepackage{amsmath}
\newcommand{\bm}[1]{\mbox{\boldmath${#1}$}}

\newcommand{\x}{\bm{x}}

\usepackage{Latex-document}

\title{\bf Fast Algorithms for Optimal Control,\vskip -2mm Anisotropic
Front Propagation and\vskip -2mm Multiple Arrivals \vskip 6mm}

\author{J. A. Sethian\vspace*{-0.5cm}\thanks{
Department of Mathematics,
University of California, Berkeley, California, USA.
E-mail: sethian@math.berkeley.edu}}
\date{\vspace{-8mm}}

\begin{document}
\maketitle

\thispagestyle{first} \setcounter{page}{735}

\begin{abstract}\vskip 3mm
We review some recent work in fast, efficient and accurate
methods to compute viscosity solutions and
non-viscosity solutions to static Hamilton-Jacobi equations which
arise in optimal control, anisotropic front propagation, and
multiple arrivals in wave propagation. For viscosity solutions,
the class of algorithms are known as ``Ordered Upwind Methods'', and
rely on a systematic ordering inherent in the characteristic flow of
information. For non-viscosity multiple arrivals, the techniques hinge
on a static boundary value phase-space formulation which again can be
solved through a systematic ordering.

\vskip 4.5mm

\noindent {\bf 2000 Mathematics Subject Classification:} 65N06, 65M06, 86A22, 49L25.

\noindent {\bf Keywords and Phrases:} Hamilton-Jacobi equations, Fast marching methods, Ordered upwind methods.
\end{abstract}

\vskip 12mm

\section{Introduction} \label{section 1}\setzero
\vskip-5mm \hspace{5mm}

This paper reviews recent work on algorithms for
static Hamilton-Jacobi equations of the form $H(Du,x)=0$; the
solution $u$ depends on $x \in R^n$, and boundary conditions
are supplied on a subset of $R^n$. These equations arise in
such areas as wave propagation, optimal control, anisotropic
front propagation, medical imaging, optics, and robotic navigation.
We develop algorithms to solve these equations
remarkably quickly, with the same optimal efficiency as classic algorithms for
shortest paths on discrete weighted networks, but extended to
continuous Hamilton-Jacobi equations.

The algorithms, which rely on
a close examination of the flow of information inherent in static
Hamilton-Jacobi equations, are robust, unconditionally stable without
time step restriction, and efficient.
They are ``One-pass'' schemes, in that
the solution is computed at $N$ grid points in $O(N \log N)$ steps.

\subsection{Viscosity vs. non-viscosity solutions}
\vskip-5mm \hspace{5mm}

\begin{figure}[hhhh]
\centerline{
$
\begin{array}{cc}
\psfig{file=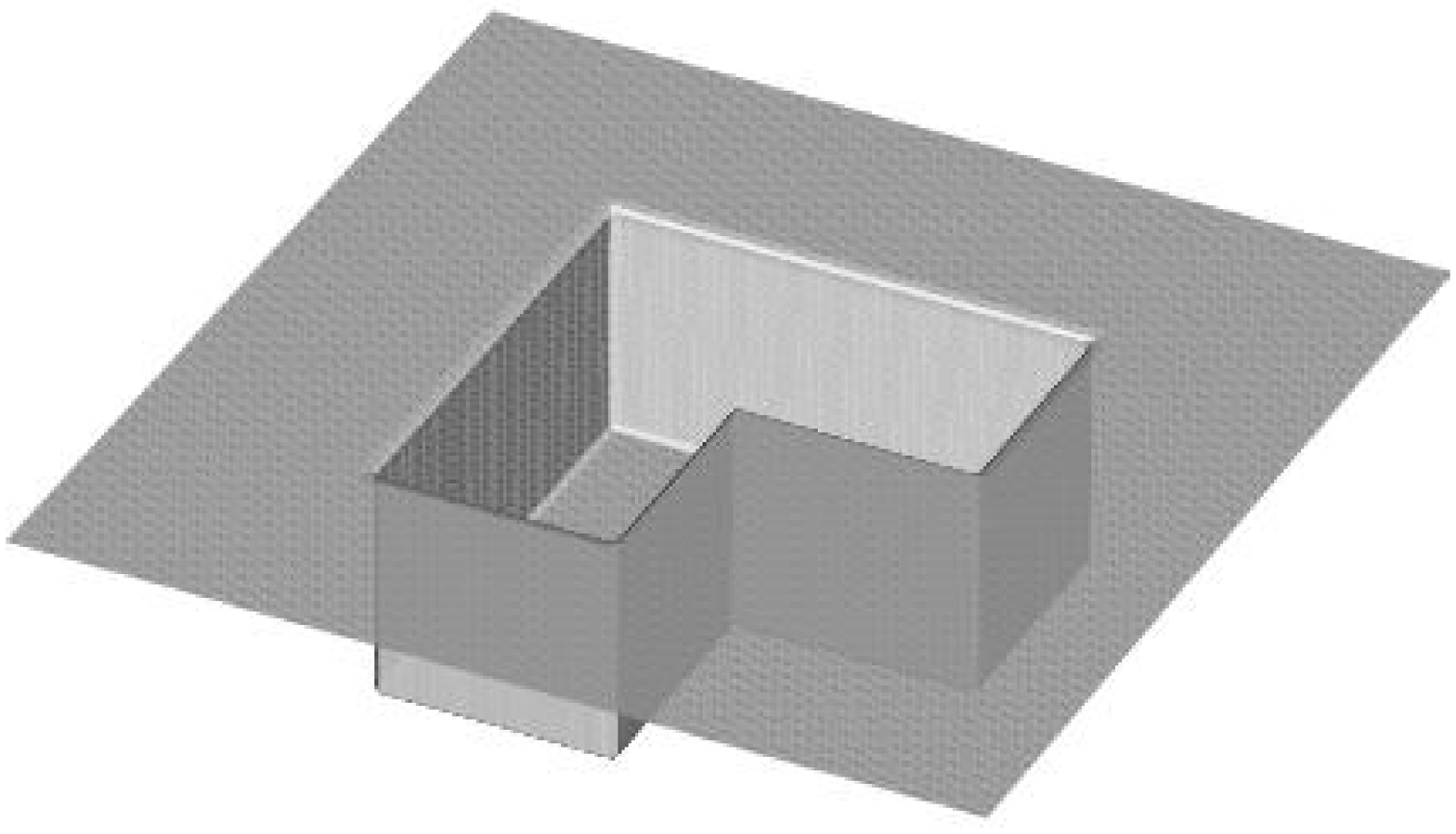,height=2in} \quad & \quad
\psfig{file=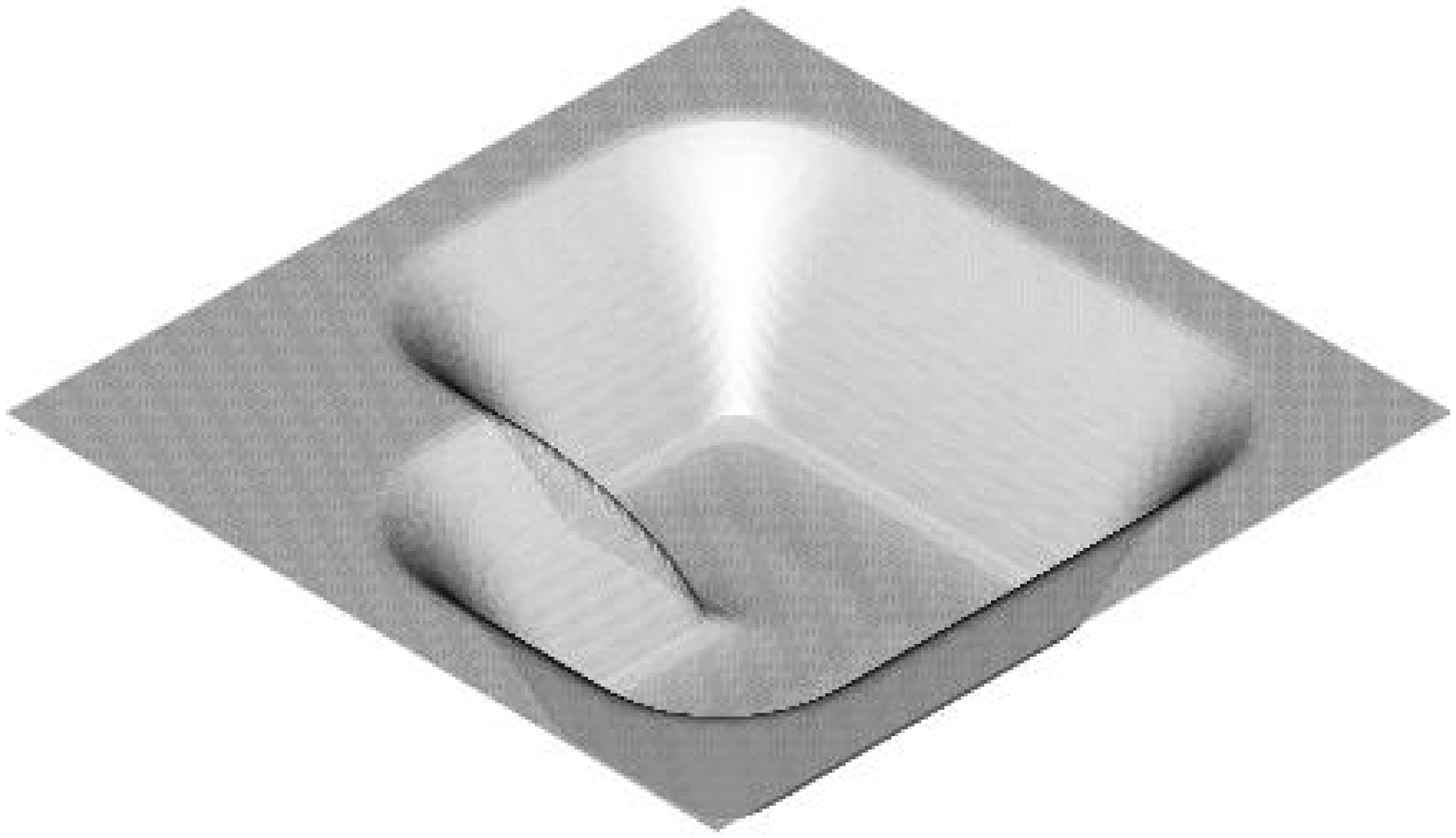,height=2in}\\
\end{array}
$ }

\centerline{
$ \mbox{\rm{Ion etching in anisotropic front
propagation}}
$
}

\centerline{ $ \mbox{\rm{Fig. \ref{figure_ion}a}} $ }

\centerline{
$
\begin{array}{cc}
\psfig{file=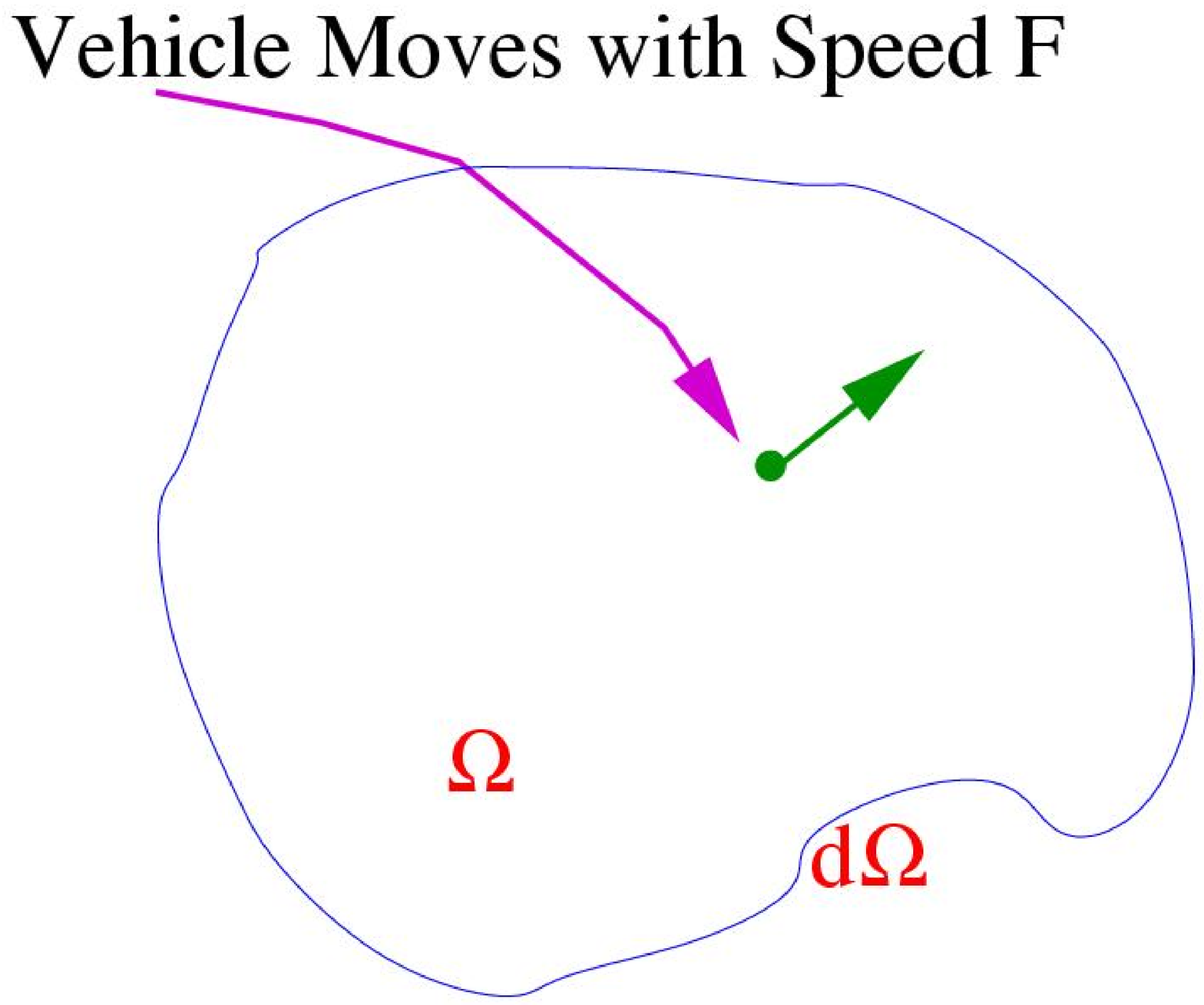,height=1.4in,angle=-90.} \quad & \quad
\psfig{file=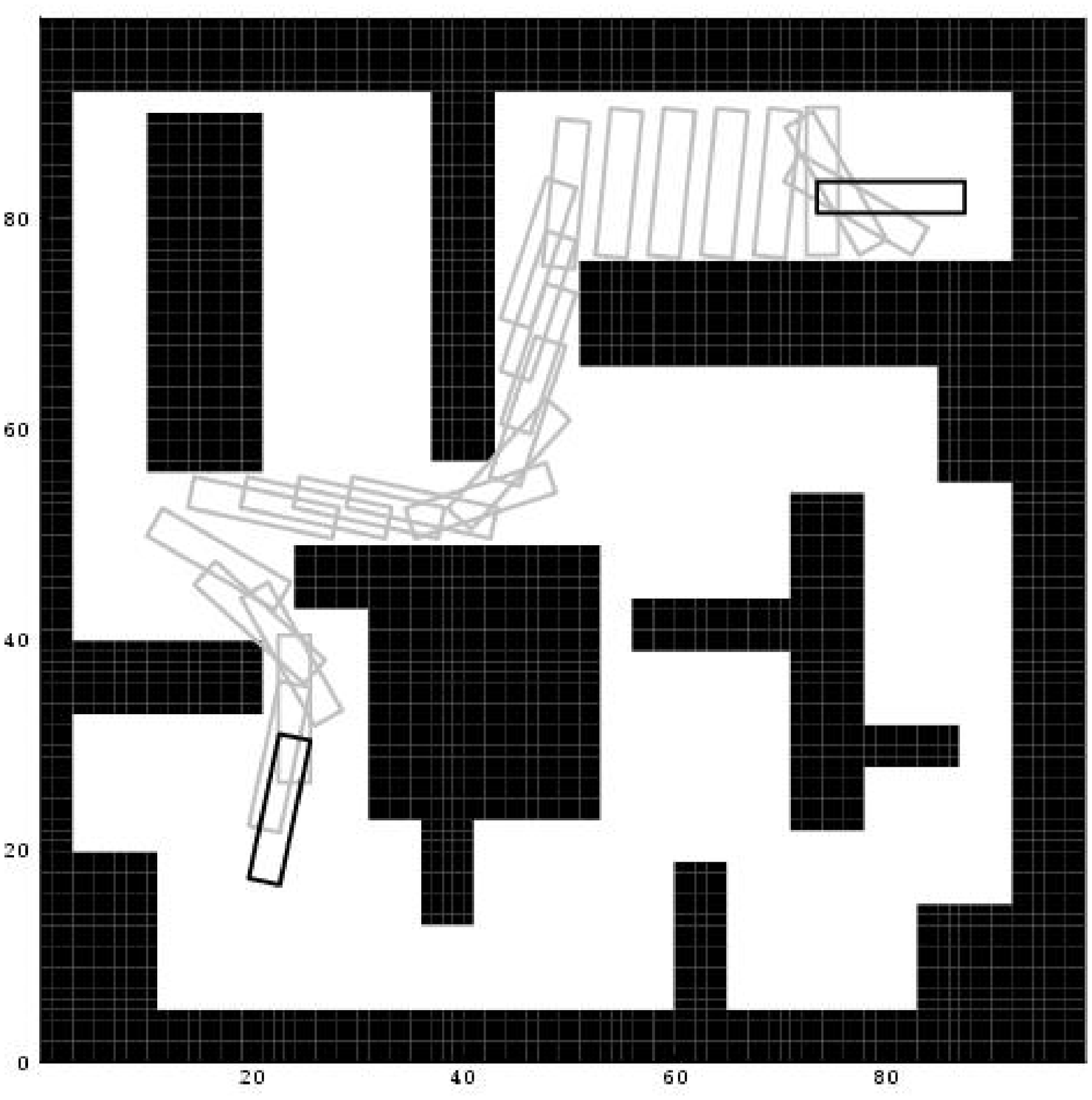,height=2.1in}\\
\end{array}
$
}

\centerline{ $ \mbox{\rm{Optimal control}} $ }

\centerline{ $ \mbox{\rm{Fig. \ref{figure_ion}b}} $ }
\caption{Viscosity solutions to static HJ equations}
\label{figure_ion}
\end{figure}

What is meant by a {\it{solution}} to $H(Du,x)=0$? Viscosity
solutions \cite{CranLion} provide a unique and well-posed
formulation which is linked to the unique viscosity limit of the
associated smoothed equation; these are first arrivals in the
propagation of information. Fig. \ref{figure_ion}a shows an
example from semiconductor manufacturing in which a beam whose
strength is angle-dependent is used to anistropically etch away a
metal surface. Fig. \ref{figure_ion}b shows an optimal control
problem to find the shortest exit path for a vehicle with position
and direction-dependent speed.

\begin{figure}[hhhh]
\centerline{ $ \psfig{file=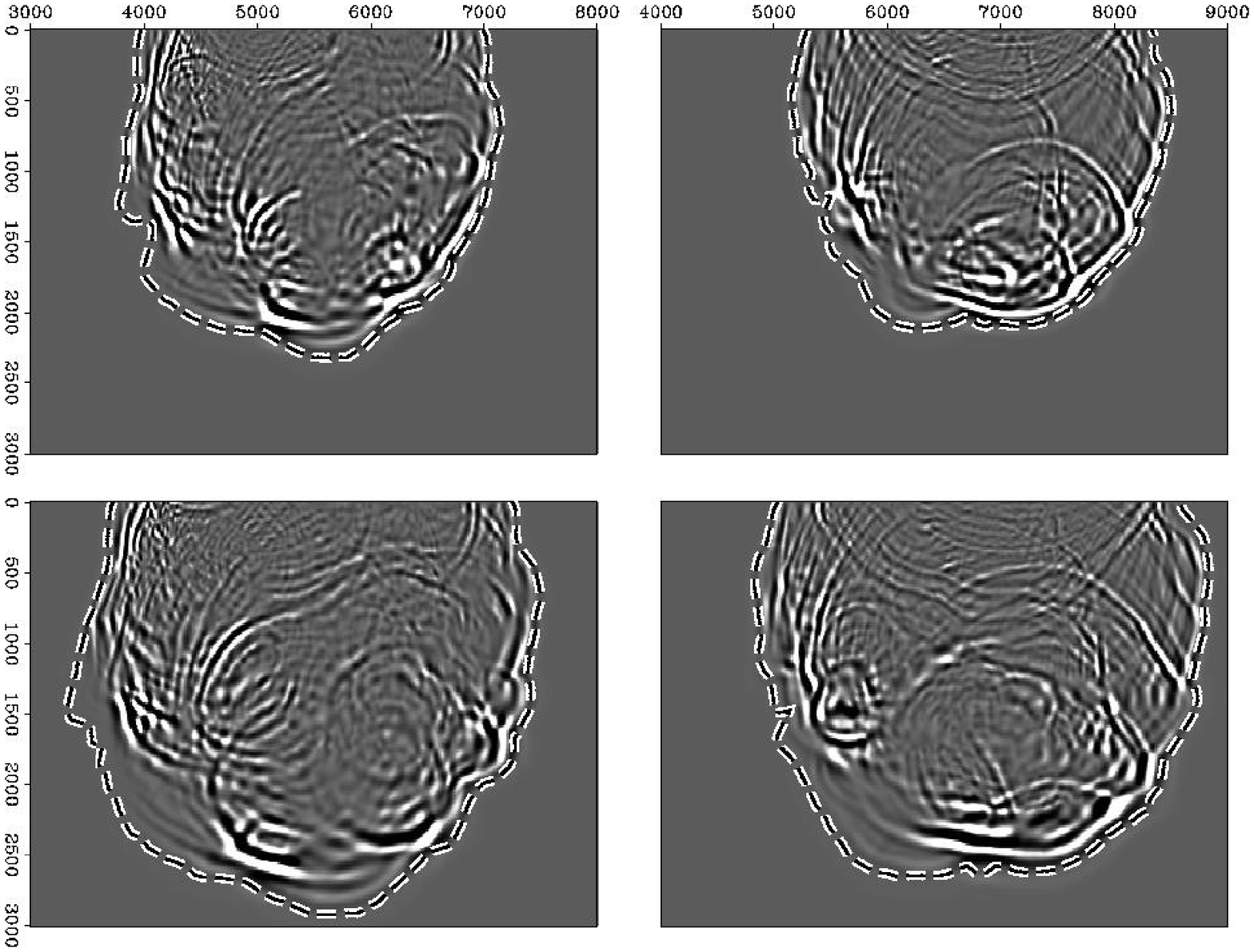,height=3.00in} $ } \centerline{ $ \mbox{\rm{Fig. \ref{figure_wave}b}} $ }
\caption{Non-viscosity solutions} \label{figure_wave}
\end{figure}

The above are viscosity solutions. However, there are many cases
in which later arrivals, or ``non-viscosity'' solutions, are
desirable. Fig. \ref{figure_wave}a shows the propagation of a wave
inwards from a square boundary; the evolving front passes through
itself and later arrivals form cusps and swallowtails as they
move; Fig. \ref{figure_wave}b shows multiple arrivals in
geophysical wave propagation.

\centerline{\psfig{file=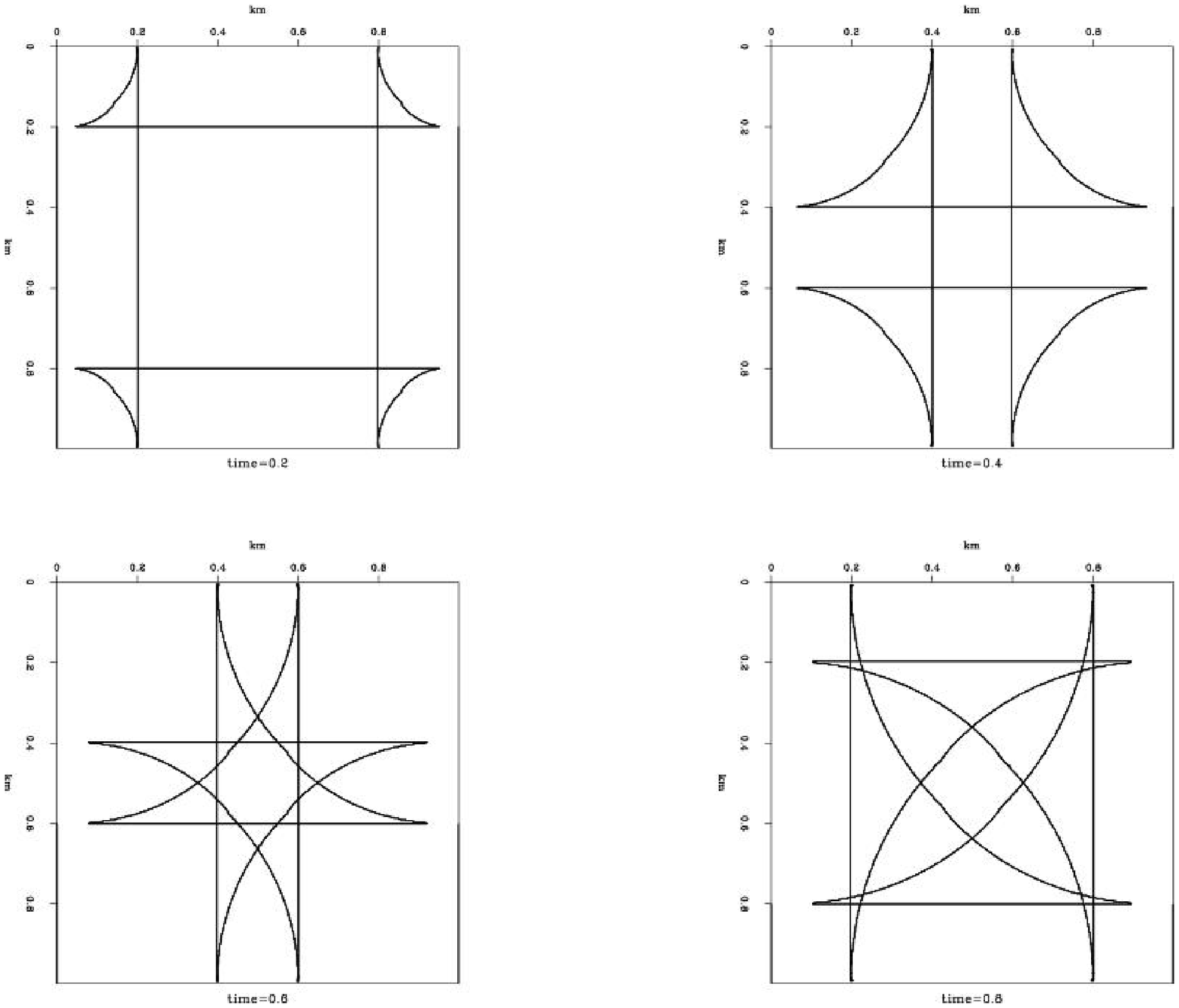,height=3.00in}}

\centerline{Fig. \ref{figure_wave}a}

Our goal is to create efficient algorithms which allow us to compute both
types of solutions. In the case of
viscosity solutions, algorithms are provided by the class of ``Ordered Upwind
Methods'' developed by Sethian and Vladimirsky in
\cite{SethVlad2, SethVlad3};
these methods work in physical space and construct the solution in
a ``One-pass'' manner through
a careful adherence to a causality inherent in the characteristic
flow of the information. In the case of non-viscosity solution,
algorithms are provided by the time-independent phase-space
formulation developed by Fomel and Sethian \cite{FomelSethian}, which relies
on conversion of multiple arrivals into an Eulerian static boundary value
problem, which can also be solved very efficiently in a
``One-pass'' manner which avoids all iteration through a careful ordering
procedure.
The remainder of this paper is devoted to describing these two classes of
algorithms and providing a few computational results.

\section{Fast methods for viscosity solutions}
\label{section 2}\setzero \vskip-5mm \hspace{5mm}

We first discuss ``Ordered Upwind Methods'' introduced
in \cite{SethVlad2} for computing viscosity solutions.

\subsection{Discrete control: Dijkstra's method}
\vskip-5mm \hspace{5mm}

Consider a discrete optimal trajectory problem on a network. Given
a network and a cost associated with each node, the global optimal
trajectory is the most efficient path from a starting point to
some exit set in the domain. Dijkstra's classic algorithm
\cite{Diks} computes the minimal cost of reaching any node on a
network in $O( N \log N)$ operations. Since the cost can depend on
both the particular node, and the particular link, Dijkstra's
method applies to both {\it{isotropic}} and {\it{anisotropic}}
control problems. The distinction is minor for discrete problems,
but significant for continuous problems. Dijkstra's method is a
``one-pass'' algorithm; each point on the network is updated a
constant number of times to produce the solution. This efficiency
comes from a careful analysis of the direction of information
propagation and stems from the optimality principle.

We briefly summarize Dijsktra's method, since the flow logic
will be important in explaining our Ordered Upwind Methods.
For simplicity, imagine a rectangular grid of size $h$,
where the cost $C_{ij} > 0$ is given for passing through each grid point
$x_{ij} = (ih, jh)$. Given a starting point, the minimal total
cost $U_{ij}$ of arriving at
the node $x_{ij}$ can be written
in terms of the minimal total cost of arriving at its neighbors:
\begin{equation}
U_{ij} = \min \left( U_{i-1,j},U_{i+1,j},U_{i,j-1},U_{i,j+1} \right)
+ C_{ij}.
\label{equation_dijkstra}
\end{equation}

To find the minimal total cost, Dijkstra's method divides mesh
points into three classes: $Far$ (no information about the correct
value of $U$ is known), $Accepted$ (the correct value of $U$ has
been computed), and $Considered$ (adjacent to $Accepted$). The
algorithm proceeds by moving the smallest $Considered$ value into
the $Accepted$ set, moving its $Far$ neighbors into the
$Considered$ set, and recomputing all {\it Considered} neighbors
according to formula \ref{equation_dijkstra}. This algorithm has
the computational complexity of $O( N \log(N))$; the factor of
$\log(N)$ reflects the necessity of maintaining a sorted list of
the $Considered$ values $U_i$ to determine the next $Accepted$
mesh point. Efficient implementation can be obtained using
heap-sort data structures.

\subsection{Continuous control: ordered upwind methods}
\vskip-5mm \hspace{5mm}

Consider now the problem of continuous optimal control; here, the
goal is to find the optimal path from a starting position to an
exit set. Dijkstra's method does not converge to the continuous
solution as the mesh becomes finer and finer, since (see
\cite{SethBook2}) it produces the solution to the partial
differential equation $ \max(|u_x|,|u_y|) = h*C, $ where $h$ is
the grid size. As $h$ goes to zero, this does not converge to the
solution of the continuous Eikonal problem given by $ |u_x^2 +
u_y^2 | ^{1/2} = C. $ Thus, Dikstra's method cannot be used to
obtain a solution to the continuous problem.

\subsubsection{Ordered upwind solvers for continuous isotropic control}
\vskip-5mm \hspace{5mm}

In the case of isotropic cost functions in which the cost depends only
on position and not on direction, two recent
algorithms, first Tsitsiklis's Method \cite{Tsitsiklis}
and then Sethian's Fast Marching Method \cite{SethFast}
have been introduced to solve the problems with the same computational
complexity as Dijkstra's method. Both methods exploit
information about the flow of information to obtain this efficiency;
the causality allows one to build the solution in increasing order, which
yields the Dijkstra-like nature of the solutions.
Both algorithms result from a key feature of Eikonal equations, namely
that their characteristic lines
coincide with the gradient lines of the viscosity solution $u(x)$; this
allows the construction of one-pass algorithms.
Tsitsiklis' algorithm evolved from studying isotropic min-time optimal
trajectory problems, and involves solving a minimization problem to
update the solution.
Sethian's Fast Marching Method evolved from studying
isotropic front propagation problems,
and involves an upwind finite difference formulation to update the
solution.
Each method starts with a particular (and different) coupled
discretization and each shows that the resulting system can be
decoupled through a causality property.
We refer the reader to these references for details on ordered
upwind methods for Eikonal equations, as well as \cite{SethVlad3}
for a detailed discussion about the similarities and differences between
the two techniques.

\subsubsection{Ordered upwind solvers for continuous anisotropic general
optimal control}
\vskip-5mm \hspace{5mm}

Consider now the full continuous optimal control problem, in which the cost
function depends on both position and direction.
In \cite{SethVlad2, SethVlad3}, Sethian and Vladimirsky
built and developed single-pass ``Ordered Upwind Methods'' for any
continuous optimal control problem. They showed how to
to produce the solution $U_i$ by recalculating each $U_i$ at most $r$ times,
where $r$ depends only
the equation and the mesh structure,
but not upon the number of mesh points.

Building one-pass Dijkstra-like methods for
general optimal control is
considerably more challenging than
it is for the Eikonal case, since characteristics no longer coincide
with gradient lines of the
viscosity solution.
Thus, characteristics and gradient lines may in fact lie in different
simplexes.
This is precisely why both Sethian's Fast Marching Method and Tsitsiklis'
Algorithm cannot be directly applied in the anisotropic (non-Eikonal)
case:
it is no longer possible to de-couple the system by computing/accepting
the mesh points in the ascending order.

The key idea introduced in
\cite{SethVlad2, SethVlad3} is to use the location anisotropy of the
cost function to limit of the number of points on the accepted front
that must be examined in the update of each Considered point.
Consider the
anisotropic min-time optimal trajectory problems, in which
the speed of motion depends not only on position but also on direction.
The value function $u$ for such problems is the viscosity solution
of the static Hamilton-Jacobi-Bellman equation
\begin{equation}
        \begin{array}{ll}
        \max_{a \in S_{1}} \left\{ (\nabla u(x) \cdot (-a)) f(a, x) \right\} = 1
,
        & x \in \Omega,\\
        u(x) = q(x), & x \in \partial \Omega.
        \end{array}
\label{eq:general2}
\end{equation}
In this formulation, $a$ is the unit vector determining the direction of
motion,
$f(a, x)$ is the speed of motion in the direction $a$ starting from the point
$x \in \Omega$, and $q(x)$ is the time-penalty for exiting the domain at
the point $x \in \partial \Omega$.  The maximizer $a$ corresponds to
the characteristic direction for the point $x$.
If $f$ does not depend on $a$, Eqn. \ref{eq:general2} reduces to the
Eikonal equation, see \cite{Bellman}.

Now, define the
anisotropy ratio $F_1/F_2$, where
$
0 < F_1 \leq f(a, x) \leq F_2 < \infty
$.
In \cite{SethVlad3}, two key lemmas were proved:
\begin{itemize}
\item{
{\bf Lemma 1.} {\it Consider the characteristic passing through
$\bar{x} \in \Omega$ and level curve $u(x) = C$, where $q_{max} <
C < u(\bar{x})$. The characteristic intersects that level set at
some point $\tilde{x}$.  If $\bar{x}$ is distance $d$ away from
the level set then $ \|\tilde{x} - \bar{x}\| \leq d
\frac{F_2}{F_1}. $ }
}
\item{
{\bf Lemma 2.} {\it Consider an unstructured mesh $X$ of diameter
$h$ on $\Omega$. Consider a simple closed curve $\Gamma$ lying
inside $\Omega$ with the property that for any point $x$ on
$\Gamma$, there exists a mesh point $y$ inside $\Gamma$ such that
$\|x-y\| < h$. Suppose the mesh point $\bar{x_i}$ has the smallest
value $u(\bar{x_i})$ of all of the mesh points inside the curve.
If the characteristic passing through $\bar{x_i}$ intersects that
curve at some point $\tilde{x_i}$ then $ \|\tilde{x_i} -
\bar{x_i}\| \leq h \frac{F_2}{F_1}. $ }
}
\end{itemize}

Thus, one may use the anisotropy ratio to exclude a large fraction
of points on the Accepted Front in the update of any Considered
Point; the size of this excluded subset depends on the anisotropy
ratio. Building on these results, a fast, Dijkstra-like method was
constructed. As before, three of mesh points classes are used. The
{\it Accepted Front} is defined as a set of $Accepted$ mesh
points, which are adjacent to some not-yet-accepted mesh points.
Define the set $AF$ of the line segments $x_j  x_k$, where $x_j$
and $x_k$ are adjacent mesh points on the $AcceptedFront$, such
that there exists a $Considered$ mesh point $x_i$ adjacent to both
$x_j$ and $x_k$. For each $Considered$ mesh point $x_i$ one
defines the part of $AF$ ``relevant to $x_i$'':
$$NF(x_i) = \left\{ (x_j, x_k) \in AF \; |
\exists \tilde{x} \mbox{ on } (x_j, x_k) \mbox{ s.t. }
\|\tilde{x} - x_i\| \leq h \frac{F_2}{F_1} \right\}.$$
We will further assume that
some consistent upwinding update formula is available:  if the
characteristic
for $x_i$ lies in the simplex $x_i x_j x_k$
then $U_i = K(U_j, U_k, x_i, x_j, x_k)$.
For the sake of notational simplicity we will refer to this value as
$K_{j,k}$.

\begin{enumerate}
\item {Start with all mesh points in $Far$ ($U_i = \infty$).}
\item {Move the boundary mesh points ($x_i \in \delta \Omega$)
to $Accepted$  ($U_i = q(x_i)$).}
\item {Move all the mesh points $x_i$ adjacent to the boundary into
$Considered$ and evaluate the tentative value of
$U_i = \min_{(x_j, x_k) NF(x_i)} K_{j,k}$.}
\item {Find the mesh point $x_r$ with the smallest value of $U$ among
all the $Considered$.}
\item {Move $x_r$ to $Accepted$ and update the {\it Accepted Front}.}
\item {Move the $Far$ mesh points adjacent to $x_r$ into $Considered$.}
\item {Recompute the value for all the $Considered$ $x_i$ within the
distance $h \frac{F_2}{F_1}$ from
$x_r$.  If less than the previous
tentative value for $x_i$ then update $U_i$.}
\item {If $Considered$ is not empty then go to 4).}
\end{enumerate}

\subsubsection{Analysis and results}
\vskip-5mm \hspace{5mm}

This is a ``single-pass" algorithm since the maximum number of
times each mesh point can be re-evaluated is bounded by the number
of mesh points in the $h \frac{F_2}{F_1}$ neighborhood of that
point; the method formally has the computational complexity of
$O((\frac{F_2}{F_1})^2 M \log(M))$. Convergence of the method to
the viscosity solution is proved in \cite{SethVlad3}, and depends
on the upwinding update formula $U_i = K(U_j, U_k, x_i, x_j,
x_k)$.

\begin{figure}[hhhh]
\centerline
{
$
\begin{array}{cc}
\psfig{file=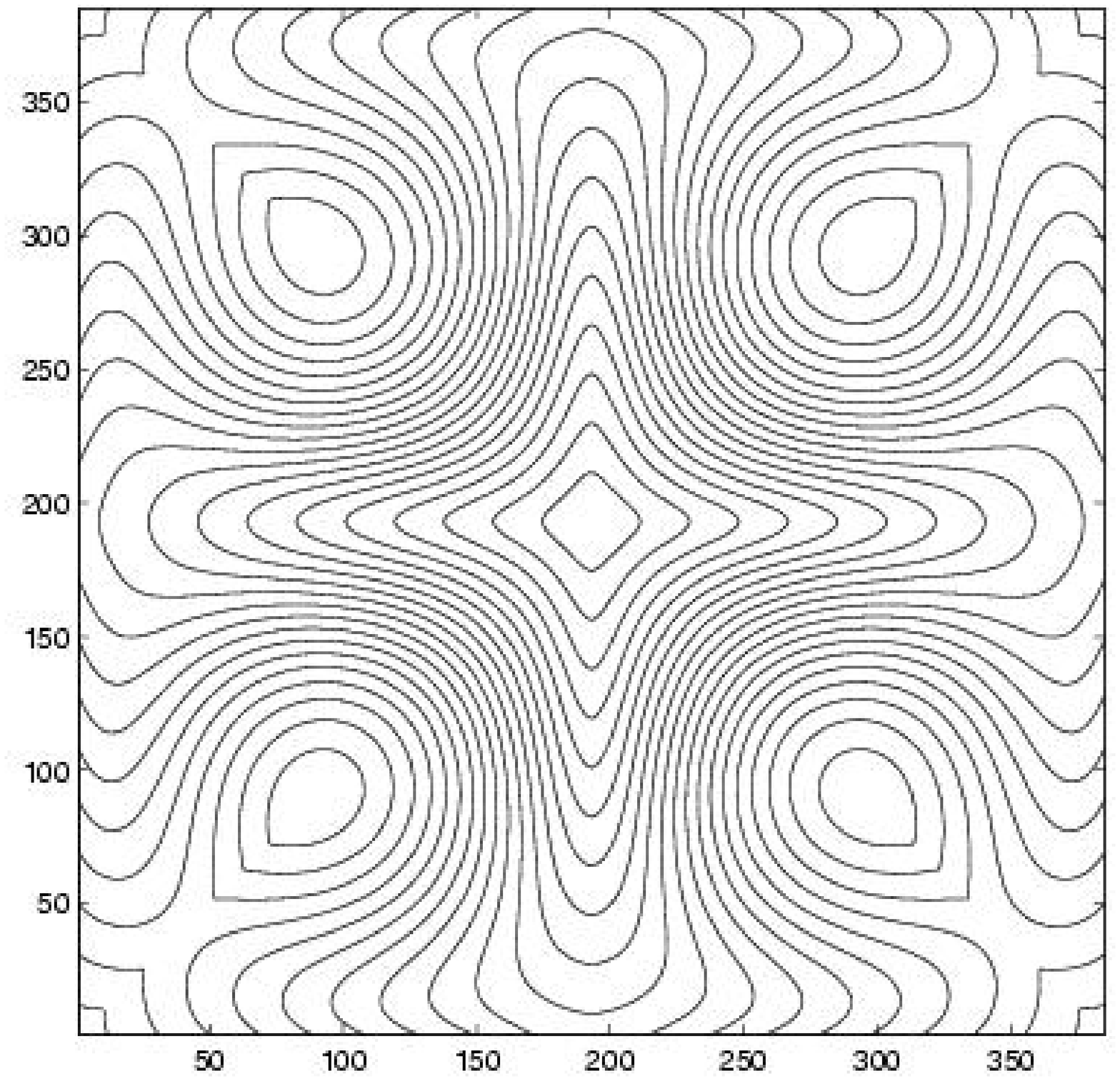,height=7cm} & \hspace{8mm}
\psfig{file=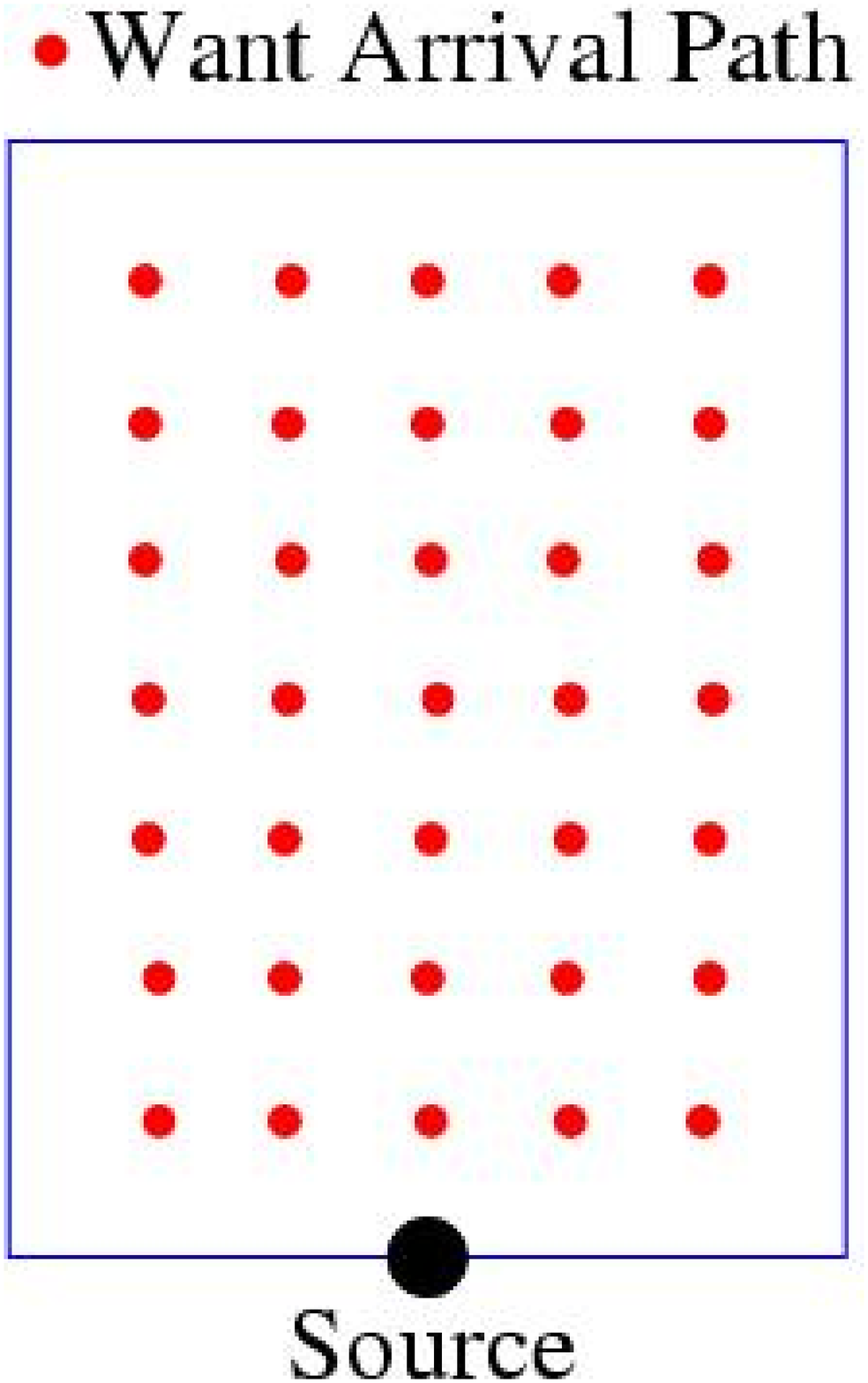,height=2.625in}\\
\end{array}
$ } \caption{Left: Anistropic front propagation \hspace{.1in}
Right: Arrival paths} \label{figure_front_ani}
\end{figure}

As an example, taken from \cite{SethVlad2},
we compute the geodesic distance on the manifold
$g(x,y) = .9 \sin(2 \pi x) \sin(2 \pi y)$
from the origin.
This can be shown to be equivalent to solving the
static Hamilton-Jacobi equation
\begin{equation}
\| \nabla u(\x) \|
F \left( \x, \frac{\nabla u(\x)}{\|\nabla u(\x)\|} \right) = 1,
\end{equation}
with speed function $F$ given by
$
F (x, y, \omega) = \sqrt{
\frac{1 + g_y^2 \cos^2(\omega) + g_x^2 \sin^2(\omega) - g_x g_y \sin(2
\omega)}
{1 + g_x^2 + g_y^2}}
$
where $\omega$ is the angle between $\nabla u(x, y)$ and the positive
direction of the $x$-axis.
The anisotropy is
substantial, since the dependence of $F$ upon $\omega$ can be pronounced
when $\nabla g$ is relatively large.
Equidistant contours are shown on the left in Figure \ref{figure_front_ani}.

\section{Fast methods for multiple arrivals} \label{section 3}\setzero
\vskip-5mm \hspace{5mm}

\subsection{Computing multiple arrivals}
\vskip-5mm \hspace{5mm}

We now consider the problem of multiple arrivals. As an example, consider the two-dimensional Eikonal equation
\begin{equation}
|\nabla u| F(x,y)=0
\label{equation-Eikonal}
\end{equation}
with $F(x,y)$ given.
We imagine a computational domain and a source point, as
shown on the right in Figure \ref{figure_front_ani}.
Suppose the goal is to determine the arrival time and path to each point
in the interior from the source point. Here, we are interested not
only in the first arrival, but all later arrivals as well.

One popular approach to computing multiple arrivals is to work in
phase space, in which the dimensionality of the problem is
increased from physical space to include the derivative of the solution as
well.
There are two approaches to computing these multiple arrivals through
a phase space formulation.
One is a Lagrangian (ray tracing) approach, in which
the phase space characteristic equations
are integrated, often from a source point, resulting in a Lagrangian
structure which fans out over the domain.
Difficulties can occur in either in low ray density
zones where there are very few rays or near caustics where rays
cross.
The other is an Eulerian description of the
problem, in either the physical domain or phase space.
In recent years, this has led to many fascinating and clever
Eulerian PDE-based approaches to computing
multiple arrivals, see, for example, \cite{Symes,Stein,Ruuth,Enqg, Benamou}.
We note that the regularity of the phase space has been utilized previously in
theoretical studies on the asymptotic wave propagation
\cite{Maslov}.
The above phase space approaches to solving for multiple arrivals
have two characteristics in common:

\begin{itemize}
\item{A phase space formulation increases the dimensionality
of the problem. In two physical dimensions, the phase space formulation
requires three dimensions; in three physical dimensions, the phase
space formulation is in five dimensions.}
\item{Given particular sources, the problem
is solved with those source location(s) as initial data. Different
sources requires re-solving the entire problem.}
\item{The problem is cast as an initial
value partial differential equation, and is evolved in time. Time step
considerations in regions of high velocity play a role in the stability of
the underlying scheme.}
\end{itemize}

\subsection{A boundary value formulation}
\vskip-5mm \hspace{5mm}

Fomel and Sethian \cite{FomelSethian} take a different approach. A
set of time-independent ``Escape Equations'' are derived, each of
which is an Eulerian boundary value partial differential equation
in phase space.  Together, they give the exit time, location and
derivative of all possible trajectories starting from all possible
interior points. Thus, the particular choice of sources is reduced
to post-processing. The computational speed depends on whether one
wants to obtain results for all possible boundary conditions, or
in fact only for a particular subset of possibilities.

\subsubsection{Liouville formulation}
\vskip-5mm \hspace{5mm}

Briefly (see \cite{FomelSethian} for details) begin with the
static Hamilton-Jacobi equation
\begin{equation}
H( x, \nabla u ) = 0,
\label{eq:hamilton-jacobi}
\end{equation}
and write the well-known characteristic equations
in phase space $(x, p)$, where $p$ corresponds to $\nabla u$
(see, for example, \cite{Evans}).
The characteristics must obey
\begin{equation}
  \label{eq:ham-jac}
  \frac{d x}{d \sigma}  =  \nabla_p H; \hspace{.5in}
  \frac{d p}{d \sigma}  =  - \nabla_x H .
\end{equation}
Differentiating the function $u(x(\sigma))$, we obtain an additional
equation for transporting the function $u$ along the characteristics:
\begin{equation}
  \label{eq:dudx}
  \frac{d u}{d \sigma} = \nabla u \cdot \frac{d x}{d \sigma} = p \cdot
   \nabla_p H .
\end{equation}
Eqns. \ref{eq:ham-jac},\ref{eq:dudx} can be initialized
at $\sigma=0$: $x(0)=x_0$, $p(0)=p_0$, $u(0)=0$.

One can now convert the phase space approach into a set of Liouville equations.
To simplify notation, we denote the phase-space vector
$(x,p)$, by $y$, the
right-hand side of system given in Eqn. \ref{eq:ham-jac} by
vector function $R(y)$, and the right-hand side of
Eqn. \ref{eq:dudx} by the function $r(y)$. In this notation,
the Hamilton-Jacobi system is
\begin{equation}
\frac{\partial y(y_0,\sigma)}{\partial \sigma}  =  R(y); \hspace{.2in}
 \frac{\partial u(y_0,\sigma)}{\partial \sigma}  =  r(y),
\label{eq:ode}
\end{equation}
and is initialized at $\sigma=0$ as $y=y_0$ and $u=0$.
This system
satisfies
\begin{equation}
  \label{eq:liou}
  \frac{\partial y(y_0,\sigma)}{\partial \sigma} = \nabla_0 y \, R(y_0)\;,
\end{equation}
and the transported function~$u$ satisfies the analogous equation
\begin{equation}
  \label{eq:liou2}
  \frac{\partial u(y_0,\sigma)}{\partial \sigma} =
\nabla_0 u \, R(y_0) + r(y_0)\;,
\end{equation}
\noindent
where $\nabla_0$ denotes the gradient with respect to $y_0$.
These are the Liouville equations.

\subsubsection{Formulation of escape equations}
\vskip-5mm \hspace{5mm}

The key idea in \cite{FomelSethian} is as follows. Assume a closed
boundary $\partial \mathcal{D}$ in the $y$ space that is crossed
by every characteristic trajectory originating in $y_0 \in
\mathcal{D}$. This defines for every $y_0$ the function $\sigma =
\widehat{\sigma}(y_0)$ of the first crossing of the corresponding
characteristic with $\partial \mathcal{D}$. Now introduce a
differentiable function $\Gamma(y)$ that identifies the boundary,
that is, $\Gamma(y) = 0$. In particular, we then have that
$\Gamma\left(y(y_0,\widehat{\sigma}(y_0)\right) = 0$. One can then
differentiate with respect to the initial condition $y_0$ to
obtain an escape equation for the parameter $\widehat{\sigma}$.
Similarly, one can derive escape equations for the position and
value, yielding the full set of

\begin{eqnarray}
\nonumber {\bf{\underline{Escape\ Equations}}} \hspace{.2in}
1 + \nabla_0 \widehat{\sigma} \cdot R(y_0) = 0\\
\nonumber
\nabla_0 \widehat{y}\,R(y_0) =0 \\
\nonumber
\nabla_0 \widehat{u} \cdot R(y_0) + r(y_0) = 0\\
\label{eq:escape_equations}
\end{eqnarray}

\subsection{Fast solution of escape equations}
\vskip-5mm \hspace{5mm}

Summarizing, rather than compute in physical space, we derive
boundary value Escape equations in phase space $y=(x,p)$. All time
step considerations are avoided, and one can compute all the
arrivals from all possible sources simultaneously. This Eulerian
formulation means that the entire domain is covered, even quiet
slow zones.

Finally, and most importantly,
a constructive, ``One-pass'' algorithm, similar to the one presented
for viscosity solutions, can be designed. Exit time, position, and derivative
at the boundary form boundary conditions.
We can then systematically march the solution inwards
in phase space from the boundary, constructing the solution through an
ordering sequence based on the characteristics that ensures
computational phase space mesh points need not be revisited more than once.

Consider a square boundary as an example, and suppose we wish to find the
time $\widehat{u} (x,z,\theta)$ at which a ray leaving the initial
point $(x,z)$ inside
the square, initially moving in direction $\theta$, hits the boundary.
We assume that the slowness field $n(x,z)$ is given.
First, note that the set $\widehat{u}(x,z,\theta) = T$,
drawn in $x,z,\theta$ space, gives the set of all initial positions
and directions which reach the boundary of the square at time $T$.
By the uniqueness of characteristics,
the set of all points parameterized by $T$ and given by
$
\widehat{U} (T)  = \{ x, z,\theta \hspace{.10in} | \hspace{.10in}
\widehat{u}(x,z,\theta) = T \}
$
sweep out the solution space.
Figure \ref{figure.surface}a shows the solution surfaces
$\widehat{u}(x,z,\theta)$ for the collapsing square.

Details on the exact algorithm are given in \cite{FomelSethian}. As demonstration (see \cite{FomelSethian}), in
Figure \ref{figure.surface}b, the top pair shows all the arrivals starting from a source at the center of the top
wall, together with the slowness field on the right (darker is slower). The bottom pair shows the first arrival
and on the amplitude of the displayed arrival (the lighter the tone, the more amplitude).

\label{lastpage}

\end{document}